\documentstyle[12pt]{article}
\begin{document}
\begin{center}
{\bf SOLUTION OF GENERALIZED FRACTIONAL REACTION-\\
  DIFFUSION EQUATIONS}\\[1cm]
R.K. SAXENA\\
Department of Mathematics and Statistics, Jai Narain Vyas University\\
Jodhpur-342005, India\\[0.5cm]
A.M. MATHAI\\
Department of Mathematics and Statistics, McGill University\\
Montreal, Canada H3A 2K6\\[0.5cm]
H.J. HAUBOLD\\
Office for Outer Space Affairs, United Nations\\
P.O.Box 500, A-1400 Vienna, Austria
\end{center}
\bigskip
\noindent
{\bf Abstract.} This paper deals with the investigation of a closed form  solution of a generalized  fractional reaction-diffusion equation. The solution of the proposed problem is developed in a compact form in terms of the H-function by the application of direct and inverse Laplace and Fourier transforms. Fractional order moments and the asymptotic expansion of the solution are also obtained.

\section{Introduction} 
It is a known fact that reaction-diffusion models play a very
 important role in pattern formation in biology, chemistry. and physics, see Wilhelmsson and Lazzaro (2001) and Frank (2005). These systems indicate that diffusion can produce the spontaneous formation of spatio-temporal patterns. For details, one can refer to the work of Nicolis and Prigogine (1977) and Haken (2004). A general model for reaction-diffusion systems is investigated by Henry and Wearne (2000, 2202) and Henry et al. (2005). 

     The simplest reaction-diffusion models are of the form 
\begin{equation}
\frac{\partial N}{\partial t}= d\frac{\partial^2 N}{\partial x^2}+ F(N), N = N (x,t), 
\end{equation} 
where $d$ is the diffusion constant and $F(N)$ is a nonlinear function representing  reaction kinetics. It is interesting to observe that for $F(N)=\gamma V(1-N)$, eq.(1) reduces to Fisher-Kolmogorov  equation  and if we set
$F(N)=\gamma N(1-N^2)$, it gives rise to the real Ginsburg-Landau equation. Del-Castillo-Negrete, Carreras and Lynch (2002) studied the front propagation and segregation in a system of reaction-diffusion equations with cross-diffusion. Recently, del-Castillo-Negrete et al (2003) discussed the dynamics in reaction-diffusion systems with non-Gaussian diffusion caused by asymmetric L\'evy flights and solved the following model
\begin{equation}
\frac{\partial N}{\partial t} = \eta D^\alpha_x N+F(N), N=N(x,t),
\end{equation} 
with $F=0$. 

     In this paper, we present a solution of a more general model of fractional  reaction-diffusion systems (2) in which $\frac{\partial N}{\partial t}$  has been replaced by the Riemann-Liouville fractional derivative $_0D_t^\beta, \beta>0$. The results derived are of a more  general nature and than those investigated  earlier by many authors, notably by Jespersen, Metzler, and  Fogedby (1999) for  anomalous diffusion  and del-Castillo-Negrete  et al. (2004) for the reaction-diffusion systems with L\'evy flights,  and fractional diffusion equation by  Kilbas et al. (2005). The solution has been developed  in terms of the H-function  in  a  compact  and elegant form with the help of Laplace and Fourier transforms and their inverses. Most of the results obtained are in a form suitable for numerical computation. The present study is in continuation of our earlier works, Haubold (1998), Haubold and Mathai (2000), and Saxena, Mathai, and Haubold (2002, 2004a, 2004b, 2004c, 2005).

\section{Mathematical Preliminaries} 
A generalization of the Mittag-Leffler function (Mittag-Leffler, 1903, 1905),
\begin{equation}
E_\alpha(z)=\sum^\infty_{n=0}\frac{z^n}{\Gamma(n\alpha+1)}, \alpha \in C, Re(\alpha)>0,
\end{equation}
was introduced by Wiman (1905)  in the generalized form 
\begin{equation}
E_{\alpha, \beta}(z)=\sum^\infty_{n=0}\frac{z^n}{\Gamma(n\alpha+\beta)},(\alpha,\beta \in c, Re(\alpha)>0, Re(\beta)>0.
\end{equation}
The main results of these functions are available in the handbook of Erd\'{e}lyi, Magnus, Oberhettinger, and Tricomi (1955, Section 18.1) and monographs  by Dzherbashyan (1966, 1993).  
     The H-function is defined by means of a Mellin-Barnes type integral in the following manner (Mathai and Saxena, 1978),
\begin{eqnarray}
H_{p,q}^{m,n}(z)&=&H^{m,n}_{p,q}\left[z\left|^{(a_p,A_p)}_{(b_q,B_q)}\right.\right]\nonumber\\
&=&H^{m,n}_{p,q}[z\left|^{(a_1,A_1)\ldots,(a_p,A_p)}_{(b_1,B_1)\ldots,(b_q,B_q)}\right.]=\frac{1}{2\pi i}\int_\Omega\Theta(\xi)z^{-\xi}d\xi,
\end{eqnarray}                
where,
 $i=(-1)^{1/2}$, 
\begin{equation}
\theta(\xi)=\frac{[\Pi^m_{j=1}\Gamma(b_j+B_j\xi)][\Pi^n_{j=1}\Gamma(1-a_j-A_j\xi)]}{[\Pi^q_{j=m+1}\Gamma(1-b_j-B_j\xi][\Pi^p_{j=n+1}\Gamma(a_j+A_j\xi)]},
\end{equation}                                    
and an empty product is always interpreted as unity; $m, n, p, q \in N_0$ with  
$0\leq n\leq p, 1\leq m\leq q, A_j, B_j \in R_+, a_i, b_j \in R$ or $ C(i=1,\ldots,p; j=1,\ldots,q)$ such that
\begin{equation}
A_i(b_j+k)\neq B_j(a_i-l-1),k,l \in N_0; i=1,\ldots, n; j=1,\ldots, m,
\end{equation}
where we employ the usual notations: $N_0=(0,1,2,\ldots,); R=(-\infty, \infty), R_+=(0,\infty)$ 
and C being the complex number field.  The contour $\Omega$ is either $L_{-\infty}, L_{+\infty}$ or $L_{i\gamma\infty}$. The following are the definitions of these contours.

(i)$\;\;\Omega=L_{-\infty}$ is a left loop situated in a horizontal strip starting at the point $-\infty+i \varphi_1$ 
and terminating at the point $-\infty+i\varphi_2$ with $-\infty<\varphi_1<\varphi_2<+\infty$;\\ 

(ii)$\;\;\Omega=L_{+\infty}$  is a right loop situated in a horizontal strip starting at the point $+\infty+i\varphi_1$ and terminating at the point $+\infty+i\varphi_2$ with $-\infty<\varphi_1<\varphi_2<+\infty$.\\

(iii)$\;\;\Omega=L_{i\gamma\infty}$  is a contour starting at the point $\gamma-i\infty$  and terminating at the point  $\gamma+i\infty$, where $\gamma\in R=(-\infty,+\infty).$

    A detailed and comprehensive account of the H-function is available from the monograph  by Mathai and Saxena (1978), Prudnikov et al.  (1989), and   Kilbas and Saigo (2004). The relation connecting the $_p\Psi_q(z)$  and the H-function is given for the first time in the monograph  by Mathai and Saxena 
(1978, p.11, Eq.1.7.8) as
\begin{equation}
_p\Psi_q[^{(a_1,A_1),\ldots,(a_p, A_p)}_{(b_1,B_1),\ldots,(B_q,B_q)}\left| z\right.]=H^{1,p}_{p,q+1}[-z\left|^{(1-a_1,A_1),\ldots,(1-a_p,A_p)}_{(0,1),(1-b_1,B_1),\ldots,(1-b_q,b_q)}\right.],
\end{equation}
where $_p\Psi_q(z)$ is  Wright`s generalized hypergeometric function, (Wright,
(1935, 1940); also see Erd\'{e}lyi, Magnus, Oberhettinger,  and Tricomi (1953, Section 4.1), defined by means of the series representation in the form
\begin{equation}
_p\Psi_q(z)=_p\Psi_q[^{(a_p,A_p)}_{(b_q,B_q)}\left|z\right.]=\sum^\infty_{r=0}\frac{[\Pi^p_{j=1}\Gamma(a_j+A_jr)]z^r}{[\Pi^q_{j=1}\Gamma(b_j+B_jr)](r)!},
\end{equation}
where $z\in C, \;a_i,b_j\in C,A_i,B_j\in R_+ ;\;i=1,\ldots,p; j=1,\ldots, q,$\\
$\sum^q_{j=1}b_j-\sum^p_{j=1}A_j>-1; C$ being the complex number field. This function includes many special functions besides the Mittag-Leffler function defined by  the equations  (3) and (4). It is interesting to observe that for $A_i=B_j=1\; \forall\; i\;\mbox{and}\; j$, (9) reduces to a generalized hypergeometric function $_pF_q(z)$. Thus

\begin{equation}
_p\Psi_q[^{(a_p,1)}_{(b_q,1)}\left|z\right.]=\frac{\Pi^p_{j=1}\Gamma(a_j)}{\Pi^q_{j=1}\Gamma(b_j)}\;_pF_q(a_1,\ldots,a_p; b_1,\ldots,b_q; z),
\end{equation}                
where $a_j\neq -\nu\;(j=1,\ldots,p$ and $\nu=0,1,2,\ldots); p<q\; (\mbox{or}\; p=q,|z|<1).$

Prior to (9), Wright (1933) introduced a special case of (9) in the form 
$$\Phi(a,b;z)=\;_0\Psi_1[^-_{(b,a)}\left|z\right.]=\sum^\infty_{r=0}\frac{1}{\Gamma(ar+b)}\frac{z^r}{(r)!},$$        
which widely occurs in problems of fractional  diffusion. It has been shown by Saxena, Mathai, and Haubold (2004b) that 
\begin{eqnarray}
E_{\alpha,\beta}(z)&=&_1\Psi_1[^{(1,1)}_{(\beta,\alpha)}\left|z\right.]\\
&=& H^{1,1}_{1,2}[-z\left|^{(0,1)}_{(0,1)(1-\beta,\alpha)}\right.].
\end{eqnarray}                                                  
If we further take $\beta=1$ in (11) and (12), we find that
\begin{eqnarray}
E_{\alpha,1}(z)&=&E_\alpha(z)=\;_1\Psi_1[^{(1,1)}_{(1,\alpha)}|z]\\ 
&=&H^{1,1}_{1,2}[-z\left|^{(0,1)}_{(0,1),(0,\alpha)}\right.],
\end{eqnarray}
where $Re(\alpha)>0, \alpha\in C.$
The definitions of the well-known Laplace and Fourier transforms  of a function $N(x,t)$ and their inverses  are described  below.

The Laplace transform of a function $N(x,t)$ with respect to $t$  is defined by 
\begin{equation}
L\left\{N(x,t)\right\}=\int^
\infty_0 e^{-st}N(x,t)dt,\;\; t>0,x\in R, 
\end{equation} 
where $Re(s) > 0$, and its inverse transform with respect to  $s$  is given by
\begin{equation}
L^{-1}\left\{N(x,t)\right\}=\frac{1}{2\pi i}\int^{\gamma+i \infty}_{\gamma-i\infty}e^{st}N(x,s)ds,
\end{equation}
$\gamma$                                                 
 being a  fixed real  number. 

     The Fourier transform of a function $N(x,t)$ with respect to $x$ is defined by 
\begin{equation}
F\left\{N(x,t)\right\}=\int^\infty_{-\infty} e^{ikx}N(k,t)dk.
\end{equation}
      
     The inverse Fourier transform with respect to $k$ is given by the eq.
\begin{equation}
F^{-1}\left\{N(x,t)\right\}=\frac{1}{2\pi}\int^\infty_{-\infty}e^{-ikx}N(k,t)dk.
\end{equation}
                                                  
     The space of functions for which the transforms defined by (15) and (17) exist  is denoted by $LF=L(R_+)\times F(R).$
                                                               
     In view of the results of Saxena,  Mathai, and Haubold (2004a, p.49), also see
Prudnikov, Brychkov, and Marichev (1989, p.355, Eq.2.25.3.2), the cosine transform of the H-function is given by 
\begin{eqnarray}
&\int^m_0&t^{\rho-1} cos(kt)H^{m,n}_{p,q}[at^\mu\left|^{(a_p, a_p)}{(b_q,B_q)}\right.]dt\\\nonumber
&=&\frac{\pi}{k^\rho}H^{n+1,m}_{q+1,p+2}\left[\frac{k^\mu}{a}\left|^{(1-b_q, B_q),(\frac{1+\rho}{2}, \frac{\mu}{2})}_{(\rho,\mu),(1-a_p,A_p),(\frac{1+\rho}{2},\frac{\mu}{2})}\right.\right],
\end{eqnarray}
where $Re[\rho+\mu^{\mbox{min}}_{1\leq j\leq m}(\frac{b_j}{B_j})]>1, |arg\; a|<\frac{1}{2}\pi \theta; \theta=\sum^n_{j=1}A_j-\sum^p_{j=n+1}A_j+\sum^m_{j=1}B_j-\sum^q_{j=m+1}B_j>0.$

     The Riemann-Liouville fractional integral of order $\nu$ is defined by 
 (Miller and Ross, 1993, p. 45)
\begin{equation}
_0D_t^{-\nu} N(x,t)=\frac{1}{\Gamma(\nu)}\int_0^t(t-u)^{\nu-1}N(x,u)du,
\end{equation}
where $Re(\nu)>0.$

     Following Samko, Kilbas, and Marichev (1990, p. 37), we define the fractional derivative  of order $\alpha>0$     in the form 
\begin{equation}
_0D_t^\alpha N(x,t)=\frac{1}{\Gamma(n-\alpha)}\frac{d^n}{dt^n}\int_0^t\frac{N(x,u)du}{(t-u)^{\alpha-n+1}};\;\ t>0, (n=[\alpha]+1,
\end{equation}
where $[\alpha]$   means the integral part of the number $\alpha$.
 
A comprehensive account of the various applications of fractional calculus in physics is available from the monograph  of Frank (2005). 
 
     From Erd\'{e}lyi et al. (1954a and 1954b, p.182), we have 
\begin{equation}
L\left\{_0D_t^{-\nu} N(x,t)\right\}=s^{-\nu}F(x,s),
\end{equation}                                                     
where $F(x,s)$ is the Laplace transform with respect to $t$ of $N(x,t), Re(s)>0$ and $Re(\nu)>0$. 
                                                                 
      The Laplace transform of the fractional derivative, defined by (21), is given by  Oldham and Spanier (1974, p. 134, Eq.(8.1.3))
\begin{equation}
L\left\{_0D_t^\alpha N(x,t)\right\}=s^\alpha N(x,s)-\sum^n_{r=1}s^{r-1}\;_0D_t^{\alpha-r}N(x,t)|_{t=0}, n-1<\alpha\leq n.
\end{equation}
In certain boundary value problems the following fractional derivative of order $\alpha$ is introduced by Caputo (1969) in the form
\begin{eqnarray}
_0D_t^\alpha f(x,t)&=&\frac{1}{\Gamma(m-\alpha)}\int_0^t\frac{f^{(m)}(x,\tau)d\tau}{(t-\tau)^{\alpha+1-m}}\\
&m-1&<\alpha\leq m, Re(\alpha)>0,\;m\in N,\nonumber\\
&=&\frac{\partial^m f(x,t)}{\partial t^m}, \mbox{if}\;\alpha=m,
\end{eqnarray}
where $\frac{d^m}{dt^m}f$ is the $m^{th}$ derivative of  order $m$  of  the function $f(x,t)$ with respect to $t$. 
     The Laplace transform of this derivative is given by Podlubny (1999) in the form  
\begin{equation}
L\left\{_0D_t^\alpha f(x,t);s\right\}= s^\alpha F(x,s)-\sum^{m-1}_{r=0} s^{\alpha-r-1}f^{(r)}(x,0+), m-1<\alpha \leq m.
\end{equation}
The above formula is  useful in deriving the solution of differential and integral equations of fractional order governing certain physical problems of reaction and diffusion . 
     
We also need the Weyl fractional derivative, defined by
\begin{equation}
_{-\infty}D^\mu_x f(x,t)=\frac{1}{\Gamma(n-\mu)}\frac{d^n}{dx^n}\int_{-\infty}^x\frac{f(u,t)du}{(x-u)^{\mu-n+1}},
\end{equation}
where $x\in R, \mu>0, n=[\mu]+1, [\mu]$ being the  integer  part of $\mu>0$
(Samko et al, 1990, Section 24.2).\\
Its Fourier transform is given by (Metzler and Klafter, 2000, p.59, A.11)
\begin{equation}
F\left\{_{-\infty}\;D_x^\mu f(x,t)\right\}= (ik)^\mu \Psi(k,t), \mu>0,
\end{equation}
where $\Psi(k,t)$   is the Fourier transform of $f(x,t)$ with respect to the variable
$x$ of  $f(x,t)$. Following the convention initiated by Compte (1996), we suppress the imaginary  unit  in Fourier space by adopting the slightly modified form of above result in our investigations(Metzler and Klafter, 2000, p.59, A.12)
\begin{equation}
F\left\{_{-\infty}D_x^\mu f(x,t)\right\}=-|k|^\mu\Psi(k,t)
\end{equation}
instead of (28). Finally we also need the following property of the H-function 
(Mathai and Saxena 1978)
\begin{equation}
H_{p,q}^{m,n}[x^\delta\left|^{(a_p,A_p)}_{(b_q,B_q)}\right.]=\frac{1}{\delta}H^{m,n}_{p,q}[x\left|^{(a_p,A_p/\delta)}_{(b_q,B_q/\delta)}\right.],
\end{equation}    
where  $\delta>0$.

\section{Fractional Reaction-Diffusion Equation}
In this section, we will investigate the solution of the generalized reaction-diffusion equation (31). The result is given in the form of the following\\  
{\bf Theorem.} Consider the generalized  fractional-reaction diffusion model 
\begin{equation}
_0D_t^\beta N(x,t)=\eta_{-\infty}D^\alpha_x N(x,t)+\varphi(x,t); \eta>0, t>0, x\in R, 1<\beta\leq 2,
\end{equation}
$0\leq\alpha \leq 1$ with the initial conditions 
\begin{equation}
_0D_t^{\alpha-1}N(x,0)=f(x), _0D_t^{\alpha-2}N(x,0)=g(x), x\in R, lim_{x\rightarrow\pm \infty}N(x,t)=0,
\end{equation}                                                      
where $_0D_t^{\alpha-1}N(x,0)$  means the Riemann-Liouville fractional derivative  of order $\alpha-1$  with respect to t evaluated at $t=0$. Similarly, $_0D_t^{\alpha-2}N(x,0)$  means the  Riemann-Liouville fractional derivative of order  $\alpha-2$ with respect to  t evaluated at  $t=0$. The quantity $\eta$   is a diffusion constant  and $\varphi(x,t)$  is a nonlinear function belonging to the area of reaction-diffusion. Then for the solution of (31), subject to the initial conditions (32), there holds the eq.
\begin{eqnarray}
N(x,t)&=&\frac{t^{\beta-1}}{2\pi}\int_{-\infty}^\infty\Psi(k)E_{\beta,\beta}(-\eta|k|^\alpha t^\beta)exp(-ikx)dk\\
&+&\frac{t^{\beta-2}}{2\pi}\int^\infty_{-\infty}\tilde{g} (k)E_{\beta,\beta-1}(-\eta|k|^\alpha t^\beta)exp(-ikx)dk\nonumber\\
&+&\frac{1}{2\pi}\int_0\xi^{\beta-1}\int_{-\infty}^\infty\varphi^\sim (k,t-\xi)E_{\beta,\beta}(-\eta|k|^\alpha\xi^\beta)exp(-ikx)dkd\xi,\nonumber
\end{eqnarray}
where $\sim$  indicates the Fourier transform with respect to the space variable $x$.\\ 
{\bf Proof.} If we apply the Laplace transform with respect to the time variable $t$  and use  eq. (23), eq. (31) becomes
$$s^\beta N^*(x,s)-f(x)-sg(x)=\eta_{-\infty}D_x^\alpha N^*(x,s)+\varphi^*(x,s). $$
As is customary, it is convenient to employ the symbol  $*$  to  indicate the Laplace transform with respect to the  variable  $t$.

    Now we apply the Fourier transform with respect to the space variable $x$ to the above equation and use the initial conditions and the result (29), then the above equation transforms into the form
\begin{equation}
\tilde{N}^*(k,s)=\frac{\tilde{f}(k)}{s^\beta+\eta|k|^\alpha}+\frac{\tilde{sg}(k)}{s^\beta+\eta|k|^\alpha}+\frac{\tilde{\varphi}^*(k)}{s^\beta+\eta|k|^\alpha},
\end{equation}
where, according to convention followed, $\sim$ indicates the Fourier transform with respect to  $x$.
On taking the inverse Laplace transform of (34) and using  the  result 
\begin{equation}
L^{-1}\left\{\frac{s^{\beta-1}}{a+s^\alpha};t\right\}=t^{\alpha-\beta} E_{\alpha, \alpha-\beta+1}(-at^\alpha),
\end{equation}                                                 
where $Re(s)> 0, Re(\alpha-\beta)>-1$, it is seen that 
\begin{eqnarray}
\tilde{N}(k,t)&=&\tilde{f}(k)t^{\beta-1}E_{\beta,\beta}(-\eta|k|^\alpha t^\beta)+\tilde{g}(k)t^{\beta-2}E_{\beta, \beta-1}(-\eta|k|^\alpha t^\beta)\nonumber\\
&+&\int_0\tilde{\varphi}(k,t-\xi)\xi^{\beta-1}E_{\beta,\beta}(-\eta|k|^\alpha\xi^\beta)d\xi.
\end{eqnarray}
     The required solution (33) now readily follows by taking the inverse Fourier transform of (36). Thus, we obtain
\begin{eqnarray*}
N(x,t)&=&\frac{t^{\beta-1}}{2\pi}\int_{-\infty}^\infty f^{\tilde{}}(k)E_{\beta,\beta}(-\eta|k|^\alpha t^\beta)exp(-ikx)dk\\
&+&\frac{t^{\beta-2}}{2\pi}\int_{-\infty}^\infty \tilde{g}(k)E_{\beta,\beta-1}(-\eta|k|^\alpha t^\beta)exp(-ikx)dk\nonumber\\
&+&\frac{1}{2\pi}\int_0^t\xi^{\beta-1}\int_{-\infty}\varphi^*(k,t-\xi)E_{\beta,\beta}(-\eta|k|^\alpha \xi^\beta)exp(-ikx)dkd\xi.
\end{eqnarray*}
     This completes the proof of the theorem.\\
{\bf Note 1.} It may be noted here that by virtue of the identity (12), the solution (33) can be expressed in terms of the H-function as can be seen from the solutions given in the special cases of the theorem in the next section. Further we observe  that (33) is not an explicit solution and special cases are of interest. \\

\section{Special Cases}
When $g(x) = 0$, then applying the convolution theorem of the Fourier transform to the solution (33), the theorem yields\\ 
{\bf Corollary 1.1.} The solution of fractional reaction-diffusion equation 
\begin{equation}
_0D_t^\beta N(x,t)= \eta_{-\infty} D^\alpha_x N(x,t)+ \varphi(x,t), t>0, \eta >0,
\end{equation}     
subject to the initial conditions,
\begin{equation}
_0D_t^{\alpha-1} N(x,t)|_{t=0}=f(x), _0D_t^{\alpha-2} N(x,t)|_{t=0}= 0
\end{equation}
$$\mbox{for}\; x\in R, lim_{x\rightarrow \pm\infty} N(x,t)=0,1<\beta\leq 2,$$
$0\neq\alpha\neq 1$, where $\eta$ is a diffusion, constant and $\varphi(x,t)$  is a nonlinear function belonging to the area of reaction diffusion is given by 
\begin{eqnarray}
N(x,t)&=&\int_{-\infty}^\infty G_1(x-\tau,t)f(\tau) d\tau\nonumber\\
&+& \int_0^t(t-\xi)^{\beta-1}\int_0^xG_2(x-\tau, t-\xi)\varphi(\tau, \xi)d\tau d\xi,
\end{eqnarray}     
where
\begin{eqnarray}
G_1(x,t)&=&\frac{t^{\beta-1}}{2\pi}\int_{-\infty}^\infty exp(-ikx)E_{\beta,\beta}(-\eta|k|^\alpha t^\beta)dk\nonumber\\
&=&\frac{t^{\beta-1}}{\pi\alpha}\int_0^\infty cos(kx)H^{1,1}_{1,2}[k\eta^{1/\alpha}t^{\beta/\alpha}\left|^{(0,1/\alpha)}_{(0,1/\alpha),(1-\beta,\beta/\alpha)}\right.]dk\nonumber\\
&=&\frac{t^{\beta-1}}{\alpha|x|}H^{2,1}_{3,3}[\frac{|x|}{\eta^{1/\alpha} t^{\beta/\alpha}}\left|^{(1,1/\alpha),(\beta,\beta/\alpha),(1,1/2)}_{(1,1),(1,1/\alpha),(1,1/2)}\right.],\;\;Re(\alpha)>0,
\end{eqnarray}
\begin{eqnarray}
G_2(x,t)&=&\frac{1}{2\pi}\int_{-\infty}^\infty exp(-ikx)E_{\beta, \beta}(-\eta|k|^\alpha t^\beta )dk\nonumber\\
&=&\frac{1}{\pi \alpha}\int^\infty_0 cos (kx)H^{1,1}_{1,2}[k\eta^{1/\alpha} t^{\beta/\alpha}\left|^{(0,1/\alpha)}_{(0,1/\alpha),(1-\beta,\beta/\alpha)}\right.]dk\nonumber\\
&=& \frac{1}{\alpha|x|}H^{2,1}_{3,3}[\frac{|x|}{\eta^{1/\alpha} t^{\beta/\alpha}}\left|^{(1,1/\alpha),(\beta,\beta/\alpha),(1,1/2)}_{(1,1),(1,1/\alpha),(1,1/2)}\right.],\;Re(\alpha)>0.
\end{eqnarray}
If we set $f(x)=\delta(x)$, where $\delta(x)$ is the Dirac delta function, then we arrive at \\
{\bf Corollary 1.2.} Consider the following reaction-diffusion model\\ 
\begin{equation}
\frac{d^\beta N(x,t)}{dt^\beta}=\eta\;_{-\infty}D_x^\alpha N(x,t),\;\;\eta>0,\;\;x\in R,
\end{equation}
with the initial condition  $N(x,t =0)  = \delta(x), lim_{x\rightarrow\pm\infty} N(x,t)=0,0<\beta\leq 1,$
where $\eta$   is a diffusion constant  and $\delta(x)$  is the Dirac delta function. Then  the solution of (42) under the given initial conditions is given by
\begin{equation}
N(x,t)=\frac{t^{\beta-1}}{\alpha|x|}H^{2,1}_{3,3}[\frac{|x|}{(\eta t^\beta)^{1/\alpha}}\left|^{(1,1/\alpha),(\beta,\beta/\alpha),(1,1/2)}_{(1,1),(1,1/\alpha),(1,1/2)}\right.],
\end{equation}
where $Re(\alpha)>0,\; Re(\beta) >2.$\\
When $\beta=1/2$,  the above corollary reduces to the following result.

      Consider the following reaction-diffusion model
\begin{equation}
\frac{d^{1/2}N(x,t)}{dt^{1/2}}=\eta\;_{-\infty}D_x^\alpha N(x,t), \;\eta>0,\;\;x\in R,
\end{equation}
with the initial condition $N(x,t=0)=\delta(x), \lim_{x\rightarrow \pm\infty} N(x,t) =0,$
where $\eta$  is a diffusion constant  and $\delta (x)$  is the Dirac delta function. Then  the solution of (42) under the given initial conditions is given by 
$$N(x,t)=\frac{1}{\alpha|x|t^{1/2}}H^{2,1}_{3,3}[\frac{|x|}{(\eta t^{1/2})^{1/\alpha}}\left|^{(1,1/\alpha),(1/2,1/2\alpha),(1,1/2)}_{(1,1),(1,1/\alpha),(1,1/2)}\right.],$$
where $Re(\alpha)>0, Re(\beta)>0.$\\                           
{\bf Remark 1.} The solution of  eq. (42), as given by Kilbas et al (2005) is in terms of the inverse Laplace and inverse Fourier transforms of certain functions whereas our solution of the same equation is obtained in an explicit closed form in terms of the H-function.

An interesting case is, when $\beta\rightarrow 1$, then in view of the cancellation law for the H-function (Mathai and Saxena, 1978), (43) provides the following  result  given by Jespersen et al. (1999) and  recently by del-Castillo-Negrete et al.  (2003) in an entirely  different form.

    For the solution of fractional-reaction-diffusion equation
\begin{equation}
\frac{d}{dt} N(x,t) = \eta_{-\infty} D_x^\alpha N(x,t),
\end{equation}
with initial condition 
$$N(x,t=0)=\delta(x), \;\;lim_{x\rightarrow \pm\infty} N(x,t)=0,$$
there holds  the relation  
\begin{equation}
N(x,t)=\frac{1}{\alpha|x|}H^{1,1}_{2,2}[\frac{|x|}{\eta^{1/\alpha}t^{1/\alpha}}\left|^{(1,1/\alpha),(1,1/2)}_{(1,1),(1,1/2)}\right.],
\end{equation}
where $Re(\alpha)>0.$
In passing, it may be noted that (46) is a closed form representation of a L\'{e}vy stable law, see Metzler and Klafter (2000). It is interesting to note that as $\alpha \rightarrow 2$, the classical Gaussian solution is recovered  as 
\begin{eqnarray}
N(x,t)&=& \frac{1}{2|x|}H^{1,1}_{2,2}[\frac{|x|}{(\eta t)^{1.2}}\left.^{(1,1/2),(1,1/2)}_{(1,1),(1,1/2)}\right.\nonumber\\
&=&\frac{1}{2|x|}H^{1,0}_{1,1}[\frac{|x|}{(\eta t)^{1/2}}\left.^{(1,1/2)}_{(1,1)}\right.\nonumber\\
&=&(4\pi \eta t)^{-1/2} exp(-\frac{|x|}{4\eta t}).
\end{eqnarray}

  It is useful to study the solution (43) due to its occurrence in certain fractional   diffusion models. we will find the fractional order moments of (43) in the next section.\\
{\bf Remark 2.} Applying Fourier transform with respect to $x$ in (42), it is found that 
\begin{equation}
\frac{d^\beta}{dt^\beta}\Psi(k,t)=-\eta|k|^\alpha\Psi(k,t),\;\;0<\beta\leq 1,
\end{equation}
which is the generalized Fourier transformed diffusion equation, since for $\alpha= 2$ and  for $\beta\rightarrow 1$, it reduces to the Fourier transformed diffusion equation 
$$\frac{d\Psi(k,t)}{dt}= -\eta|k|^2\Psi(k,t),$$
being a relaxation (diffusion) equation, for a fixed wave number $k$ 
(Metzler and Klafter, 2000). Here $\Psi(k,t)$   is the Fourier transform of $N(x,t)$ with respect to $x$.\\
{\bf Remark 3.} It is interesting to observe that the method employed for deriving the solution of the equations (31) and (32) in the space $=LF=L(R_+)\times F(R)$  can also be applied  in the space $LF'=L'(R_+)\times F',$ where $F'=F'(r)$ is the space of Fourier transform of  generalized functions. As an illustration, we can choose  
$F' = S'$ or $F'=D'.$ The Fourier transforms in $S'$ and $D'$ are introduced by  Gelfand and Shilov (1964). $S'$ is the dual of the  space $S$ , which  is the space of all infinitely differentiable functions which, together with their derivatives approach zero more rapidly than any power of  $1/|x|$ as $|x|\rightarrow \infty$
(Gelfand and Shilov, 1964, p.16). $D'$ is the dual of the space $D$  which consists of all smooth functions with compact supports (Brychkov and Prudnikov, 
1989 p. 3). For  further details , the reader is referred to the monographs  by Gelfand and Shilov  (1964) and Brychkov and Prudnikov (1989),  if we replace  the Laplace and Fourier transforms in (15)  and (17) by the corresponding Laplace and Fourier transform of generalized  functions. 

\section{Fractional Order Moments} 
In this section, we will calculate the fractional order moments, defined by
\begin{equation}
<|x|^\delta>=\int_{-\infty}^\infty|x|^\delta N(x,t)dx.
\end{equation}
   
Using the definition of the Mellin transform    
\begin{equation}
M\left\{f(t);s)\right\}=\;\;\int^\infty_0t^{s-1}f(t)dt,
\end{equation}
we find from (49) that 
\begin{equation}
<|x(t)|^\delta>=\int^\infty_{-\infty}|x|^\delta N(x,t)dx
\end{equation}
\begin{equation}
<|x|^\delta(t)>=\frac{2t^{\beta-1}}{\alpha}\int^\infty_0 x^{\delta-1}H^{2,1}_{3,3}[\frac{|x|}{\eta^{1/\alpha}t^{\beta/\alpha}}\left|^{(1,1/\alpha),(1,\beta/\alpha),(1,1/2)}_{(1,1),(1,1/\alpha),(1,1/2)}\right.]dx.
\end{equation}

     Applying the Mellin transform formula for the H-function, namely 
\begin{equation}
\int^\infty_0 x^{\rho-1}H^{m,n}_{p,q}[ax\left|^{(a_p,A_p)}_{(b_q,b_q)}\right.]dx=a^{-p}\Theta(-\rho),
\end{equation}
where 
$$^{-min}_{1\leq j \leq m}Re(\frac{b_j}{B_j})<Re(\rho)<^{max}_{1\leq j \leq n} Re(\frac{1-a_j}{A_j}),\;\;|arg\; a|<\frac{1}{2}\pi\theta, \theta>0,$$
$\Theta$ is defined  in (19) and $\Theta(-\rho)$ in the definition of the H-function (5),
we see that 
\begin{equation}
<|x|^\delta(t)>=\frac{2}{\alpha}\eta^{\delta/\alpha}t^{\beta\left\{(\delta/\alpha)+1-(1/\beta)\right\}}\frac{\Gamma(-\frac{\delta}{\alpha})\Gamma(1+\delta)\Gamma(1+\frac{\delta}{\alpha})}{\Gamma(-\frac{\delta}{2})\Gamma(\beta+\frac{\beta\delta}{\alpha})\Gamma(1+\frac{\delta}{2})},
\end{equation}
whenever the gammas exist, $Re(\delta)>-1$ and $Re(\delta+\alpha) > 0.$

     Two interesting special cases of (54) are worth mentioning.
\begin{equation}
\mbox{(i) As}\;\;\delta\rightarrow 0, \mbox{then using the result}\frac{1}{\Gamma(z)}\sim z\;\;\mbox{for}\;\; z << 1,    
\end{equation}
we find that 
\begin{equation}
^{lim}_{\delta\rightarrow 0}<|x|^\delta (t)>=\frac{t^{\beta-1}}{\Gamma(\beta)}.                   
\end{equation}
(ii) When $\alpha=2, \delta=2$,  the linear time dependence,  
\begin{equation}
^{lim}_{\delta\rightarrow 2,\alpha\rightarrow 2}<|x(t)|^\delta>=\frac{2\eta t^{2\beta-1}}{\Gamma(2\beta)},
\end{equation}
of the mean squared displacement is recovered.\par

\section{Behavior of the Solution in Equation (43)}
Eq. (43) can be expressed  in terms of  the Mellin-Barnes type integral (Erd\'{e}lyi et al, 1953, chapter 1) as
\begin{equation}
N(x,t)=\frac{t^{\beta-1}}{\alpha|x|}\frac{1}{2\pi}\int_L\frac{\Gamma(\frac{s}{\alpha})\Gamma(1-s)\Gamma(1-\frac{s}{\alpha})}{\Gamma(\beta-\frac{s\beta}{\alpha})\Gamma(1-\frac{s}{2})\Gamma(s/2)}[\frac{|x|}{\eta^{1/\alpha}t^{\beta/\alpha}}]^s ds.
\end{equation}

     Let us assume that the poles of the integrand of (58) are simple. Now evaluating the sum of residues in ascending powers of $|x|$ by calculating the residues at the poles of $\Gamma(1-s)$ at the points $s=1+\nu\;(\nu=0,1,2,\ldots)$, and $\Gamma(1-\frac{s}{\alpha})$ at the points $s=1+\nu\;(\nu=0,1,2,\ldots),$ it is found that the series expansion of the general solution (43) is given by
\begin{eqnarray}      
N(x,t)&=& \frac{t^{\beta-1}}{\alpha\eta^{1/\alpha}t^{\beta/\alpha}}\sum^\infty_{\nu=0}\frac{(-1)^\nu}{(\nu)!}\frac{\Gamma[\frac{1+\mu}{\alpha}]\Gamma(1-(1+\nu)/\alpha)}{\Gamma[1-\frac{\beta(1+\nu)}{\alpha}]\Gamma[\frac{1-\nu}{2}]\Gamma[\frac{1+\nu}{2}]}\left[\frac{|x|}{\eta^{1/\alpha}t^{\beta/\alpha}}\right]^\nu\nonumber\\
&+&\frac{|x|^{\alpha-1}}{\eta t}\sum^\infty_{\nu=0}\frac{(-1)^\nu\Gamma(1-\alpha(1+\nu))}{\Gamma(-\beta\nu)\Gamma(1-\frac{\alpha}{2}(1+\nu))\Gamma(\frac{\alpha(1+\nu)}{2})}\left[\frac{|x|^\alpha}{(\eta t^\beta)^{1/\alpha}}\right]^{\alpha \nu}
\end{eqnarray}
where $0<Re(\nu)<\alpha,\left\{\frac{|x|}{\eta^{1/\alpha}t^{1/\alpha}}\right\}<1.$

From (59), we infer that
\begin{equation}
N(x,t)\sim At^{\beta-\beta/\alpha-1}+B|x|^{\alpha-1}, \mbox{as} \;x\rightarrow 0,
\end{equation}
where A and B are numerical constants .    
     Further from the series expansion (59), it can be seen that the initial behavior is given by
\begin{equation}
N(x,t)\sim\frac{\Gamma(1+\frac{1}{\alpha})\Gamma(1-\frac{1}{\alpha})}{\pi\eta^{1/\alpha}t^{1+\beta/\alpha-\beta}\Gamma(\beta-\frac{\beta}{\alpha})}\;\;\mbox{for}\;\;1<\alpha<2,
\end{equation}
 with  $\left\{\frac{|x|}{(\eta t^\beta)^{1/\alpha}}\right\}<<1.$\\
Next, if we calculate the residues at the poles of $\Gamma(s/\alpha)$   at the points\\
$s=-\alpha\nu\;(\nu=0,1,2,\ldots),$ it gives
\begin{equation}
N(x,t)=\frac{t^{\beta-1}}{|x|}\sum^\infty_{\nu=0}\frac{\Gamma(1+\alpha \nu)}{\Gamma(\beta+\beta \nu)\Gamma(1+\frac{\alpha \nu}{2})\Gamma(-\frac{\alpha \nu}{2})}\left[-\frac{\eta t^\beta}{|x|^\alpha}\right]^\nu,
\end{equation}
for $|x| > 1$. From (62), it can be readily seen that 
$$N(x,t)\sim \frac{t^{\beta-1}}{|x|}$$  
for  large $|x|.$   
In conclusion, it is observed that the solution given by (43) does not admit a probabilistic interpretation in contrast with fractional reaction-diffusion based on Caputo derivative derived by the authors. However, when  $\beta\rightarrow 1$ then  it has  a probabilistic interpretation, as shown  in corollary 1.2.\par
\bigskip
\noindent
\section*{References}
\noindent
Brychkov, Yu.A. and Prudnikov, A.P.: 1989, \emph {Integral Transforms of}
 
\emph {Generalized Functions}, Gordon and Breach, New York.\\
Caputo, M.: 1969, \emph {Elasticita e Dissipazione}, Zanichelli, Bologna.\\ 
Compte, A.: 1996, Stochastic foundations of fractional dynamics, 

\emph {Physical Review E} {\bf 53}, 4191-4193.\\
Del-Castillo-Negrete, D., Carreras, B.A., and Lynch, V.E.: 2003, Front

dynamics in reaction–diffusion systems with L\'evy flights: A fractional 

diffusion approach, \emph {arXiv: nlin.PS/0212039 v2}.\\
Del-Castillo-Negrete, D., Carreras, B.A., and Lynch, V.E.: 2002, Front

propagation and segregation in a reaction–diffusion model with

cross–diffusion, \emph {Physica D}, {\bf 168-169}, 45-60.\\
Dzherbashyan, M.M.: 1966, \emph {Integral Transforms and Representation of}

\emph {Functions in Complex Domain} (in Russian), Nauka, Moscow. \\
Dzherbashyan,M.M.: 1993, \emph {Harmonic Analysis and Boundary Value}

\emph {Problems in the Complex Domain}, Birkhaeuser–Verlag, Basel, London.\\
Erd\'{e}lyi, A., Magnus, W., Oberhettinger, F., and Tricomi, F.G.: 1953, \emph {Higher} 

{Transcendental Functions}, Vol. {\bf 1}, McGraw–Hill, New York, Toronto,

and London.\\
Erd\'{e}lyi, A., Magnus, W., Oberhettinger, F., and Tricomi, F.G.: 1954a, 

\emph {Tables of Integral Transforms}, Vol. {\bf 1}, McGraw–Hill, New York, Toronto,

and London.\\
Erd\'{e}lyi, A., Magnus, W., Oberhettinger, F., and Tricomi, F.G.: 1954b, 

\emph {Tables of Integral Transforms}, Vol. {\bf 2}, McGraw–Hill, New York, Toronto, 

and London.\\
Erd\'{e}lyi, A., Magnus, W., Oberhettinger, F., and Tricomi, F.G.: 1955, \emph {Higher}

\emph {Transcendental Functions}, Vol. {\bf 3}, McGraw–Hill, New York, Toronto, and

London.\\
Frank, T.D.: 2005, \emph {Nonlinear Fokker-Planck Equations: Fundamentals and} 

\emph {Applications}, Springer-Verlag, Berlin-Heidelberg.\\
Gelfand I.M. and Shilov, G.F.: 1964, \emph {Generalized Functions}, Vol. {\bf 1}, 

Academic Press, London.\\
Haken, H.: 2004, \emph {Synergetics: Introduction and Advanced Topics}, 

Springer-Verlag, Berlin-Heidelberg.\\
Haubold, H.J.: 1998, Wavelet analysis of the new solar neutrino capture rate

data for the Homestake experiment, \emph {Astrophysics and Space Science} {\bf 258},

201-208.\\
Haubold, H. J. and Mathai, A. M.: 2000, The fractional kinetic equation and

thermonuclear functions, \emph {Astrophysics and Space Science} {\bf 273}, 53-63.\\
Henry, B.I. and Wearne, S.L.: 2000, Fractional reaction–diffusion, \emph {Physica A}

{\bf 276}, 448-455.\\
Henry, B.I. and Wearne, S.L.: 2002, Existence of Turing instabilities in 

a two-species fractional reaction-diffusion system, 

\emph {SIAM Journal of Applied Mathematics} {\bf 62}, 870-887.\\
Henry, B.I., Langlands, T.A.M., and Wearne, S.L.: 2005, Turing pattern 

formation in fractional activator-inhibitor systems, 

\emph {Physical Review E} {\bf 72}, 026101.\\
Jespersen, S., Metzler, R., and Forgedby, H.C.: 1999,  L\'evy flights in  external

force fields: Langevin and fractional Fokker–Planck equations and their

solutions, \emph {Physical Review  E} {\bf 59}, 2736-2745.\\
Kilbas, A.A., Pierantozzi, T., Trujillo, J.J., and Vazquez, L.:2005, On 

generalized fractional evolution-diffusion equation, under publication.\\
Kilbas, A.A. and Saigo, M.: 2004, \emph {H-Transforms: Theory and Applications}, 

Chapman and Hall/CRC, New York.\\
Kulsrud, R.M.: 2005, \emph {Plasma Physics for Astrophysics}, Princeton 

University Press, Princeton and Oxford.\\
Mathai, A. M.: 1993, \emph {A Handbook of Generalized Special Functions for} 

\emph {Statistical and Physical Sciences}, Clarendon Press, Oxford.\\
Mathai, A.M. and Saxena, R.K.: 1978, \emph {The H-function with Applications in} 

\emph {Statistics and Other Disciplines}, Halsted Press [John Wiley and Sons], 

New York, London and Sydney.\\
Metzler, R. and Klafter, J.: 2000, The random walk's guide to anomalous 

diffusion: A fractional dynamics approach, \emph {Physics Reports} {\bf 339}, 1-77.\\
Metzler, R. and Klafter, J.: 2004, The restaurant at the end of the 

random walk: Recent developments in the description of anomalous 

transport by fractional dynamics, 

\emph {Journal of Physics A: Math. Gen.} {\bf 37}, R161-R208.\\
Miller, K.S. and Ross, B.: 1993, \emph {An Introduction to the Fractional Calculus} 

\emph {and Fractional Differential Equations}, John Wiley and Sons, New York.\\
Mittag-Leffler, G.M.: 1903, Sur la nouvelle function  E(x), \emph {C.R. Acad. Sci.}, 

\emph {Paris (Ser.II)} {\bf 137}, 554-558.
\\
Mittag-Leffler, G. M.: 1905, Sur la representation analytique d' une function 

branche uniforme d'une fonction, \emph {Acta Mathematica} {\bf 239}, 101-181.\\
Nicolis, G. and Prigogine, I.: 1977, \emph {Self-Organization in Nonequilibrium} 

\emph {Systems: From Dissipative Structures to Order Through Fluctuations}, 

John Wiley and Sons, New York.\\  
Oldham, K. B. and Spanier, J.: 1974, \emph {The Fractional Calculus: Theory and} 

\emph {Applications of Differentiation and Integration to Arbitrary Order}, 

Academic Press, New York, and Dover Publications, New York 2006.\\
Podlubny, I.: 1999, \emph {Fractional Differential Equations}, Academic 

Press, San Diego.\\ 
Prudnikov, A.P., Brychkov, Yu.A., and Marichev, O.I.: 1989, \emph {Integrals and} 

\emph {Series}, Vol. {\bf 3}, \emph {More Special Functions}, Gordon and Breach, New York.\\
Samko, S.G., Kilbas, A.A., and Marichev, O.I.: 1990, \emph {Fractional Integrals} 

\emph {and Derivatives: Theory and Applications}, Gordon and Breach, 

New York.\\ 
Saxena, R.K., Mathai, A.M., and Haubold, H.J.: 2002, On fractional kinetic 

equations, \emph {Astrophysics and Space Science} {\bf 282}, 281-287.\\ 
Saxena, R.K., Mathai, A.M., and Haubold, H.J.: 2004a, On generalized 

fractional kinetic equations, \emph {Physica A} {\bf 344}, 657-664.\\
Saxena, R.K., Mathai, A.M., and Haubold, H.J.: 2004b, Unified fractional 

kinetic equation and a fractional diffusion equation, \emph {Astrophysics and} 

\emph {Space Science} {\bf 290}, 299-310.\\
Saxena, R.K., Mathai, A.M., and Haubold, H.J.: 2004c, Astrophysical 

thermonuclear functions for Boltzmann-Gibbs statistics and Tsallis 

statistics, \emph {Physica A} {\bf 344}, 649-656.\\
Saxena, R.K., Mathai, A.M. and Haubold, H.J.: 2005, Fractional reaction-

diffusion equations, this volume.\\
Tsallis, C.: 2004, What should a statistical mechanics satisfy to reflect 

nature?, \emph {Physica D} {\bf 193}, 3-34.\\
Tsallis, C. and Bukman, D.J.: 1996, Anomalous diffusion in the presence of 

external forces: Exact time-dependent solutions and their 

thermostatistical basis, \emph {Physical Review E} {\bf 54}, R2197-R2200.\\
Wilhelmsson, H. and Lazzaro, E.: 2001, \emph {Reaction-Diffusion Problems in} 

\emph {the Physics of Hot Plasmas}, Institute of Physics Publishing, 

Bristol and Philadelphia.\\
Wiman, A.: 1905, Ueber den Fundamentalsatz in der Theorie der Funktionen

$E_\alpha(x)$, \emph {Acta Mathematica} {\bf 29}, 191-201.\\
Wright, E.M.: 1933, On the coefficients of power series having exponential 

singularities, \emph {Journal of the London Mathematical Society} {\bf 8}, 71-79.\\
Wright, E.M.: 1935, The asymptotic expansion of the generalized hypergeo-

metric functions, \emph {Journal of the London Mathematical Society} {\bf 10}, 

386-293.\\    
Wright, E.M.: 1940, The asymptotic expansion of the generalized hypergeo-

metric functions, \emph {Proceedings of the London Mathematical Society} {\bf 46}, 

389-408.   
\end{document}